\documentclass[journal]{IEEEtran}

\usepackage{bm}
\usepackage{wrapfig}
\usepackage{float}
\usepackage{amsmath}
\usepackage{amssymb}
\usepackage{multirow}
\usepackage{amsfonts}
\usepackage{mathrsfs}
\usepackage{graphicx}
\usepackage{epstopdf}
\usepackage{subcaption}
\usepackage{adjustbox}
\usepackage{textcomp}
\usepackage{autobreak}
\usepackage{tabularx}
\usepackage{comment}
\usepackage{caption}
\usepackage{lettrine}
\usepackage{cite}
\usepackage{setspace}
\usepackage{amsmath}
\usepackage{titlesec}
\usepackage{lipsum}

\usepackage{theoremref}

\usepackage[utf8]{inputenc}
\usepackage{tikz}
\usepackage{import}
\usetikzlibrary{shapes.geometric, arrows}
\usetikzlibrary{fit}
\tikzstyle{process2} = [rectangle, minimum width=2cm, minimum height=1cm, text centered, text width=2cm, draw=black]
\tikzstyle{process} = [rectangle, minimum width=4cm, minimum height=1cm, text centered, text width=4cm, draw=black]
\tikzstyle{arrow} = [thick,->,>=stealth]

\titlespacing*{\section}
{0pt}{3.5ex plus 1ex minus .2ex}{2.0ex plus .2ex}

\ifCLASSINFOpdf
\else
\fi
\begin{document}
%
\title{Flexibility Characterization of Sustainable Power Systems in Demand Space: A Data-Driven Inverse Optimization Approach}


\author{Mohamed~Awadalla and Fran\c{c}ois~Bouffard~\IEEEmembership{Senior Member,~IEEE}%
\thanks{This work was supported in part by IVADO, Montreal, QC, and the Natural Sciences and Engineering Research Council of Canada, Ottawa, ON.}
\thanks{M. Awadalla and F. Bouffard are with the Department of Electrical and Computer Engineering, McGill University, Montreal, QC H3A~0E9, Canada and with the Groupe d'\'{e}tudes et de recherche en analyse des d\'{e}cisions (GERAD), Montreal, QC  H3T~1J4, Canada (emails: mohamed.awadalla@mail.mcgill.ca; francois.bouffard@mcgill.ca).}
}


%
\maketitle
\begin{abstract}
The deepening of the penetration of renewable energy is challenging how power system operators cope with their associated variability and uncertainty. The inherent flexibility of dispathchable assets present in power systems, which is often ill-characterized, is essential in addressing this challenge. Several proposals for explicit flexibility characterization focus on defining a feasible region that secures operations either in generation or uncertainty spaces. The main drawback of these approaches is the difficulty in visualizing this feasibility region when there are multiple uncertain parameters. Moreover, these approaches focus on system operational constraints and often neglect the impact of inherent couplings (e.g., spatial correlation) of renewable generation and demand variability. To address these challenges, we propose a novel data-driven inverse optimization framework for flexibility characterization of power systems in the demand space along with its geometric intuition. The approach captures the spatial correlation of multi-site renewable generation and load using polyhedral uncertainty sets. Moreover, the framework projects the uncertainty on the feasibility region of power systems in the demand space, which are also called \textbf{\textit{loadability sets}}. The proposed inverse optimization scheme, recast as a linear optimization problem, is used to infer system flexibility adequacy from loadability sets. 
\end{abstract}

%
\IEEEpeerreviewmaketitle

\section*{Nomenclature}

The main symbols used in the paper are listed here. Other symbols will be defined as required.

\subsection{Sets and Indices}
\begin{IEEEdescription}[\IEEEusemathlabelsep\IEEEsetlabelwidth{$1,2,3,4$}]

\item[$\Xi_{g d}(\zeta)$]{Feasibility region of a power system in the generation-demand space.}
\item[$\Xi_{d}(\zeta)$]{Projection of the generation-demand space onto the demand space only, also called loadability set.}
\item[$\mathcal{D}(\zeta)$]{Minimal loadability set.}
\item[$\mathcal{J}$]{Set of inequalities, indexed by $j$ and of size $J$.}
\item[$\mathcal{L}$]{Set of transmission lines, indexed by $l$ and of size $L$.}
\item[$\mathcal{M}$]{Set of generating units, indexed by $m$ and of size $M$.}
\item[$\mathcal{M}_{n}$]{Set of generating units connected to bus $n$.}
\item[$\mathcal{N}$]{Set of buses, indexed by $n$ and of size $N$.}
\item[$\mathcal{T}$]{ Set of time periods, indexed by $t$ and of size $T$.}
\end{IEEEdescription}

\subsection{Variables}
\begin{IEEEdescription}[\IEEEusemathlabelsep\IEEEsetlabelwidth{$1,2,3,4$}]
\item[$d_{n}$]{Residual demand at bus $n$.}
\item[$g_{m}$]{Power output of generator $m$ that can be deployed to satisfy uncertain residual demand.}
\item[$g^{0}_{m}$]{Scheduled base point of generating unit $m$.}
\item[$\hat{g}_n$]{Total dispatchable generation output at bus $n$.}
\item[${q}_{n}$]{Net power injected at bus $n$.}
\item[$s$]{Vector of slack variables associated with strong duality constraint in the inverse problem.}
\item[$u_{m}$]{Commitment status (binary---0: if off, 1: if on) of generator $m$.}
\item[$y$]{Vector of dual variables associated with primary constraints in the forward problem.}
\item[$\zeta_m$]{Schedules of generating unit $m$ including commitment status, that is $\zeta_m$ = ($u_m$, $g_{m}^{0}+r_{m}^{\uparrow}$, $g_{m }^{0}-r_{m}^{\downarrow}$) in the operations planning horizon.}
\item[$\zeta$]{Collection of all $\zeta_m$.}
\item[$\rho$]{Goodness of fit metric for inverse optimization.}
\item[$\epsilon_{n}$]{Load shedding and renewable generation curtailment at bus $n$.}

\end{IEEEdescription}

\subsection{Parameters}
\begin{IEEEdescription}[\IEEEusemathlabelsep\IEEEsetlabelwidth{$1,2,3,4$}]
\item[$f_l^{\max}$]{Maximum flow capacity of transmission line $l$.} 
\item[$g_m^{\max}$]{Maximum power limit of generator $m$.}
\item[$g_m^{\min}$]{Minimum power limit of generator $m$.}
\item[$h_{l n}$]{Power transfer distribution factor (PTDF) of line $l$ for power injections at bus $n$}
\item[${d^{0}}$]{Vector of nominal residual demand for inverse optimization.}
\item[$\gamma$]{Value of lost load and renewable generation curtailment.} 
\end{IEEEdescription}


\section{Introduction}
\IEEEPARstart{T}{he shares} of renewable energy resources (RES) have increased significantly in the last decade and will grow at even faster paces as economies decarbonize. This large-scale, uncertain, and volatile renewable generation creates many challenges to power system operations. The ongoing transformation of the conventional power systems paradigm leads to the study of the emerging concept of power system flexibility \cite{B.Mohandes}. Power system flexibility is defined as the system’s ability to accommodate any component outage or variation in its net load (\emph{i.e.}, demand less non-dispatchable generation) to keep the system secure \cite{B.Mohandes, N.K.Dhaliwal, Y.Huo1}. The ultimate goal is to have enough flexibility to cope with the increasing RES levels so that economic and secure operation can be maintained. In the literature, several definitions and metrics have been proposed to study the need for, and the provision of power system flexibility. Broadly, quantifying power systems flexibility for power system operation with RES can be classified into two main categories: implicit and explicit methods with respect to the RES uncertainties \cite{M.Shahidehpour}. \par

The implicit approach deals with pre-defined uncertainty and solves generation scheduling problems \cite{Y.Huo1, M.Shahidehpour}. Most contributions in this category use stochastic techniques and robust approaches which model uncertain parameters using probability distributions and worst-case scenarios, respectively \cite{A.Velloso},\cite{C.wang}. However, this category of approaches only focuses on how to optimally exploit existing generation assets to deal with a given amount of uncertainty. \par

Conversely, the motivation for explicit approaches found in the literature is that future power systems are expected to have a limited flexible capacity. Furthermore, they recognize that unit commitment (UC) problem solutions may be overly costly if worst-case realizations of uncertainty have to be satisfied by dispatching expensive generators. Moreover, they acknowledge that in extreme cases it may not even be possible find feasible UC solutions \cite{C.Shao}. Thus, the practicality of classical robust UC, meant to support a wide range of operating conditions, could be severely limited under very deep penetrations of renewables. In order to address these shortfalls, the second class of approaches tackles the reverse question: how much uncertainty can be adopted using existing flexible resources? This necessitates the assessment of system flexibility for given operating conditions \cite{M.Shahidehpour}.\par

Explicit flexibility assessment is characterized mainly through the use of region-based geometrical approaches. A pioneering study proposed the concept of \emph{do-not-exceed (DNE) limits} to define the maximum variations of uncertain parameters a system can accommodate using robust optimization \cite {J.Zhao2}. These limits leverage unambiguously the utilization of renewable resources while treating uncertainty sets as a decision variable. Further enhancements of the DNE dispatch method has been proposed considering corrective topology control actions \cite{A.S.Korad}, and using historical data of wind power realizations \cite{F.Qiu}. Along similar lines, the dispatchable region concept was introduced to characterize flexibility regions explicitly. The authors of \cite {W.Wei1} characterized the dispatchable region of wind generation and revealed its geometry to be a polytope in uncertainty space. This concept was generalized in \cite {Y.Liu1} to include a full ac network model. Another form of dispatchable region was optimized using energy and reserve scheduling based on an affine redispatch policy \cite {W.Wei4}. The flexibility set determination approach proposed in \cite{M.A.Bucher} infers a polytope describing the allowed deviations from current system state in generation and tie lines spaces. Likewise, reference \cite{N.Yorino} estimated the size of a power system's feasible region in generation space under different levels of uncertainties.  \par

Similar research direction have expanded to quantify flexibility using various metrics. Zhao \emph{et  al.} \cite {J.Zhao} suggested a flexibility measure based on a system's DNE region. They defined a binary flexibility metric to check if the largest variation of uncertainty is within the admissible range or not. Another flexibility framework has quantified insufficient flexibility by power imbalance event magnitude and frequency \cite{Z.Qin}. There is a family of publications \cite{W.Wei2, C.Wang3, C.Wang4, C.Wang2} that investigated wind generation admissibility assessment using two-stage robust optimization. In these papers, wind accommodation capabilities under a given solution UC strategy is assessed using expected load shedding and wind curtailment. \par

However, most of the previous work was built on budget-constrained polyhedral uncertainty sets and neglected spatial and temporal trends in historical data \cite {W.Wei2, C.Wang3,C.Wang4}. Another drawback of budget-constrained uncertainty modeling is the combinatorial growth of vertices needed to model the uncertainty set with respect to the uncertain parameters \cite{A.Velloso}. A more accurate way of modeling uncertainty is to calculate the convex hull of spatial and temporal scenarios for wind generation admissibility assessment \cite{C.Wang2}. The empirical study in \cite{F.Golestaneh} does not recommend the use of the convex hull for wind and photovoltaics in dimensions higher than four dimensions because of inherent computational costs. \par

From a broader perspective, the majority of region-based flexibility assessment methods have assessed and visualised dispatchable regions in either generation or uncertainty spaces considering system operational constraints regardless of the uncertainty impacts of the RES output \cite{W.Wei1, Y.Liu1, W.Wei3}. Moreover, when the dimensionality of the uncertainty grows, these regions cannot be easily visualized. This paper's proposal attempts to fill this research gap by developing a comprehensive flexibility assessment framework under deep penetration of renewables, leveraging the advantages of region-based approaches and metric-based methods. In this context, the proposed method is considered a complement to energy and reserve scheduling problems \cite {Y.Huo1, A.Velloso}, and belongs to the latter type of flexibility assessments such as \cite{M.Shahidehpour}. The framework relies on a given unit commitment strategy, and consequently assesses the impact of the uncertainty of RES outputs at a given operating condition. \par

Compared to previous work, our proposal has the following features. First, it extends the notion of RES accommodation assessment proposed in \cite{J.Zhao2, W.Wei1, W.Wei2}, which maps feasibility regions from the generation-demand space to the demand space only using the loadability set approach \cite{A.Kalantari, A.A.Jahromi1, A.A.Jahromi2}. Second, it captures the spatial correlation of multiple renewable generation sites and loads using scenario-based polyhedral uncertainty sets. The framework leverages the uncertainty set and enables a reformulation of the loadability set with enhanced data-driven capabilities. In addition, a novel data-driven inverse optimization problem formulation is proposed which seeks to identify existing system flexibility for uncertainty mitigation by exploring the feasibility region of a linear programming (LP)-relaxation i.e loadability set. Finally, this assessment is meant to be quantitatively indicative of how much ``room'' exists in the bulk power system to handle residual demand uncertainty for given unit commitment solutions and network topologies. 

The main contributions of this paper are fourfold:
\begin{enumerate}
\item We develop a data-driven polyhedral uncertainty set to capture residual demand uncertainty.
\item We endogenously embodied the polyhedral uncertainty set in identifying the redundant constraints that does not shape the feasibility region in the generation-demand space. Also, we revisit the loadability set characterization and implicitly incorporate the data-driven scenario-based uncertainty set to redefine the loadability set.
\item We propose a unified framework to characterize power system flexibility explicitly and geometrically in the demand space using data-driven inverse optimization technique (DDIO).
\item We present a solution methodology for DDIO which considers uncertainty sets which may or may not intersect its feasibility region along with its geometric intuition. 

\end{enumerate}

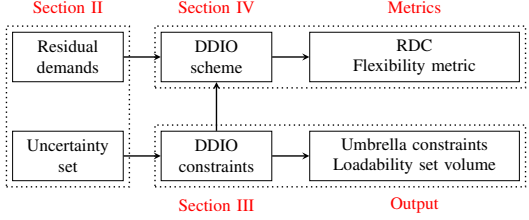
\begin{figure}[t!]
\centering
\resizebox{7cm}{!}{\begin{tikzpicture}[node distance=2cm]
\node (pro1) [process2] {Residual demands};
\node (pro2) [process2, right of=pro1,xshift=1cm] {DDIO scheme};
\node (pro3) [process, right of=pro2,xshift=2cm] {RDC\\Flexibility metric};
\node (pro4) [process2, below of=pro1] {Uncertainty
set};
\node (pro5) [process2, right of=pro4,xshift=1cm] {DDIO constraints};
\node (pro6) [process, right of=pro5,xshift=2cm] {Umbrella constraints\\ Loadability set volume};

\node (labe) [above of =pro1,yshift=-1cm,red] {Section II};
\node (labe) [above of =pro2,yshift=-1cm,red] {Section IV};
\node (labe) [above of =pro3,yshift=-1cm,red] {Metrics};
\node (labe) [below of =pro6,yshift=1cm,red] {Output};
\node (labe) [below of =pro5,yshift=1cm,red] {Section III};

\draw [arrow] (pro1) -- (pro2);
\draw [arrow] (pro2) -- (pro3);
\draw [arrow] (pro4) -- (pro5);
\draw [arrow] (pro5) -- (pro6);
\draw [arrow] (pro5) -- (pro2);

\node[draw,dotted,thick,fit=(pro5) (pro6)] {};
\node[draw,dotted,thick,fit=(pro1) (pro4)] {};
\node[draw,dotted,thick,fit=(pro2) (pro3)] {};

\end{tikzpicture}}
\caption{Overview of the method proposed in this paper}\label{fig:03}
\end{figure}

The remainder of this paper is organized as follows. Section~II describes how data-driven uncertainty sets are obtained. Section~III discusses the residual demand admissibility assessment model. Next in Section~IV, we introduce the inverse optimization technique and its application in power systems flexibility characterization. Section~V conducts two case studies to show the effectiveness of our proposed approach. Finally, Section~VI presents our conclusions. An overview of our proposal is shown in Fig.~\ref{fig:03}.


\section{Data-Driven Nodal Residual Demand Uncertainty Set}

In this section, we detail the procedure for building data-driven polyhedral uncertainty envelopes. Without loss of generality here, we focus on capturing the uncertainty of spatially-correlated residual demand. We develop a data-based polyhedral uncertainty set as per \cite{A.Velloso} by leveraging principal component analysis (PCA) \cite{Y.Huo1,A.Berizzi}. PCA is applied to two historical time series: one of observed residual demand (demand less non-dispatchable generation, also known as \emph{net load}) and a second corresponding to the previously forecasted values of those residual demands. All time series have to have the same length of $T$, and they are synchronized and evenly spaced in time (e.g., one hour interval).  \par

We denote a matrix ${W} \in {\mathbb{R}}^{T \times N}$ whose elements $w_{n t}$ are the time series of \emph{observed} residual demand at bus $n \in \mathcal{N} = \{1,...,N\}$ for each time instance $t \in \mathcal{T} = \{1,...,T\}$. Similarly, we define the matrix $\mu \in {\mathbb{R}}^{T \times N}$ whose elements $\mu_{n t}$ are the time series of \emph{forecasted} residual demand at bus $n \in \mathcal{N}$ for each past time instance $t \in \mathcal{T}$.

Using $\mu$, we obtain the centered data matrix $W_{c}$ \cite{Y.Huo1, A.Berizzi}, whose contents are the residual demand forecast errors at all nodes $n \in \mathcal{N}$ and times $t \in \mathcal{T}$
\begin{equation}
W_{c}=W-\mu
\end{equation}
Assuming residual demand forecast errors are unbiased\footnote{In the case where errors are unbiased, one would need to calculate the biases at each node and then remove them from $W_c$.}, its spatial forecast error covariance matrix $\Sigma \in {\mathbb{R}}^{N \times N}$ is approximated by
\begin{equation}
\Sigma=\frac{1}{T-1} {W}_{c}^\top {W}_{c}
\end{equation}
The next step consists in performing PCA by eigenvalue decomposition of the spatial forecast error covariance matrix.

We now let the columns of an $N \times N$ matrix, $V$, and the diagonal entries of another $N \times N$ matrix, $\Lambda$, represent, respectively, the orthonormal eigenvectors (\emph{principal components}) and the eigenvalues of the spatial forecast error covariance matrix $\Sigma$. Here, the diagonal elements of $\Lambda$ are ordered such that $\lambda_{11} \geq \lambda_{22} \geq \cdots \geq \lambda_{NN}$, while the columns of $V$ are arranged such that its $n$th column (eigenvector) is associated with the $n$th eigenvalue, $\lambda_{nn}$.

We then project the forecast errors contained in $W_c$ onto the $k = 1,\ldots,N$ eigenvectors of its covariance matrix \cite{A.Berizzi}
\begin{equation}
{Z_k}= {W}_{c} V_k
\end{equation}
The resulting $T \times 1$ vector $Z_k$ is the original forecast error data projected onto $V_k$, the $k$th principal component of $\Sigma$. The coordinates of the extrema of each principal component data projection $Z_k$ are identified
\begin{equation}
    \bar{\mathcal{S}}_k = \arg\max_{t \in \mathcal{T}} \| Z_k \|^2
\end{equation}
That is $\bar{\mathcal{S}}_k$ is the data point projected along principal component $V_k$ which is the furthest away from the origin. Keeping a conservative approach, we will assume data can range between $-\bar{\mathcal{S}}_k$ and $\bar{\mathcal{S}}_k$ along the principal component $V_k$.

Considering that typically only the first few dominant principal components (\emph{i.e.}, the ones with largest eigenvalues) are sufficient to describe accurately the forecast error uncertainty, it is common practice to limit the number of principal components to $K$ such that $K < N$. As a result, we can argue that any of the original forecast errors can be reconstructed as a convex combination of the extrema of the $K$ data projections $Z_k$\footnote{In the case where $W_c$ was biased, the forecast error biases found prior to applying PCA would also need to be added back.}
\begin{align}
    \mathcal{S}^+_k &= \bar{\mathcal{S}}_k \\
    \mathcal{S}^-_k &= - \bar{\mathcal{S}}_k
\end{align}

Inspired by the data-driven convex hull uncertainty set concept proposed previously in \cite{A.Velloso,F.Golestaneh}, we propose to define a \emph{polyhedral uncertainty set} (PUS) of historical residual demand forecast errors
\begin{align}
P(\mathcal{S}) =& \left\{\mathbf{E} \in \mathbb{R}^N \mid \mathbf{E} = \sum_{k=1}^{K} \left ( \omega^+_{k} \mathcal{S}^+_{k} + \omega^-_{k} \mathcal{S}^-_{k} \right ), \right . \nonumber \\
& \qquad \sum_{k=1}^{K} \left ( \omega^+_{k} + \omega^-_{k} \right ) = 1, \nonumber \\
& \qquad  0 \leq \omega^+_{k} \leq 1, \; 0 \leq \omega^-_{k} \leq 1,  \forall k \Bigg \}\label{eq:Pus} 
\end{align}
As can be inferred here, $P(\mathcal{S})$ is a convex and closed set whose interior contains the vast majority of the original set of forecast errors. Note, moreover, that PUS is a convex hull of the extrema of the retained $Z_k$ data projections, where we define $\mathcal{S} = \cup_k ( \mathcal{S}_k^+ \cup \mathcal{S}_k^-)$. Thereby the PUS represents the smallest convex set that contains every data point projected onto the $K$ retained principal components.

Assuming that the number of historical forecast error data points is large and reflects a good sample of a diversity of operating conditions, we argue that $P(\mathcal{S})$ is time invariant. Therefore, we posit that $P(\mathcal{S})$ shifted by a vector of nodal residual demand forecasts looking into the future (\emph{i.e.}, not historical, like $\mu$) $d^0$ is an adequate representation of possible future residual demand and its uncertainty. Otherwise said, we define the forward-looking residual demand PUS as
\begin{align}
P(\mathcal{S}, d^0) =& \left\{\mathbf{E} \in \mathbb{R}^N \mid \mathbf{E} = d^0 + \sum_{k=1}^{K} \left ( \omega^+_{k} \mathcal{S}^+_{k} + \omega^-_{k} \mathcal{S}^-_{k} \right ), \right . \nonumber \\
& \qquad \sum_{k=1}^{K} \left ( \omega^+_{k} + \omega^-_{k} \right ) = 1, \nonumber \\
& \qquad  0 \leq \omega^+_{k} \leq 1, \; 0 \leq \omega^-_{k} \leq 1,  \forall k \Bigg \}\label{eq:Pus1} 
\end{align}

\begin{figure}
\centering
\includegraphics[width=3.0in]{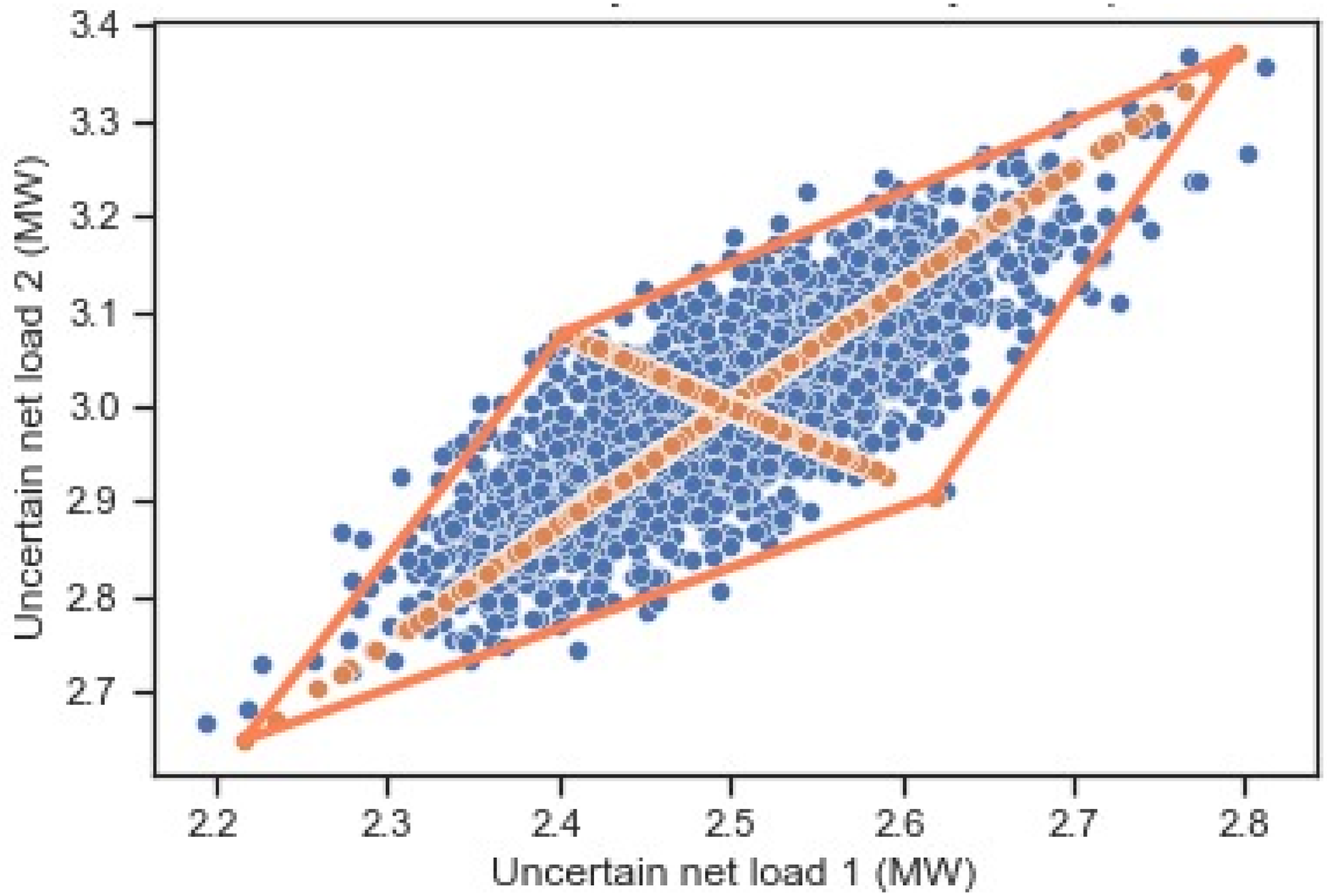}
\caption{Two-dimensional polyhedral uncertainty set.}\label{Fig_system}
\end{figure}

Fig.~\ref{Fig_system} illustrates the construction of a forward-looking PUS in a two node power system. In this example, we see the original historical residual forecast error data, represented by blue dots, and the data projected onto two of its principal components (orange dots). By inspection, the rhombus-shaped envelope, whose principal axes correspond to the principal components of the data, encapsulates the vast majority of the original data. Moreover, the vertices of this rhombus correspond to the set of points $\mathcal{S}$. Moreover, we notice here that $d^0 = \begin{bmatrix} 2.5 & 3.0 \end{bmatrix}^\top$~MW is the forecast residual demand at nodes 1 and 2, respectively.

Compared to the conventional convex hull encapsulation approach \cite{A.Velloso}, the number of edges of the PUS depends on the chosen number of retained principal components $K$, which can be at most $2N$. Conversely, the convex hull is sensitive to all possible realizations that a power system might encounter, albeit unlikely. This affects the complexity of the resulting convex hull because of data outliers, and consequently its number of vertices can grow exponentially when the dimensionality of the system uncertainty increases.

\section{Residual Demand Uncertainty Accommodation Assessment Framework}

\subsection{Introduction}

As argued in introduction, a standard two-stage robust unit commitment with a single uncertainty set \cite{A.Velloso} determines in the first-stage the generators' commitment status, energy and reserve schedules. The second-stage finds the minimum redispatch solution given the worst uncertainty realization. This computationally-heavy approach goes against existing industry practice, where stress tests are applied \emph{ex-post} as offline validation procedures \cite{A.Velloso}. 

Next, we present, starting from a first-stage unit commitment solution, how much uncertainty can be handled by a power system. For expository purposes, we use a relaxed linear optimization of the second-stage assuming first-stage variables are known, as it is usually the case in most RES admissibility assessment frameworks \cite{J.Zhao2, W.Wei1, W.Wei2, C.Wang3, C.Wang4}. 

\subsection{Problem Formulation}

A residual demand admissibility framework evaluates quantitatively how much demand less non-dispatchable renewables can be accommodated by the bulk power system given a specific UC solution without causing any operational infeasibility.

We formulate the following linear optimization problem, which optimizes residual demand admissibility as a function of the first-stage unit commitment solution, expressed by the collection of the tuples $\zeta_m = (u_{m}, g_{m}^{0}+r_{m}^{\uparrow}, g_{m}^{0}-r_{m}^{\downarrow})$ for each of the system's dispatchable generators $m \in \mathcal{M} = \{1, \ldots,M \}$. Specifically, $\zeta_m$ contains: $u_{m}$ its commitment status (on/off), $g_{m}^{0}$ its scheduled power output, and $r_{m}^{\uparrow}$ and $r_{m }^{\downarrow}$ respectively its scheduled up- and down-reserves.

\begin{equation}
f(\zeta_1,\ldots, \zeta_M) = \min_{\hat{g}, q, \epsilon, d}  \sum_{n=1}^N \gamma\left|\epsilon_{n}\right|  \label{eq:P1a}
\end{equation}
Subject to:
\begin{align}
q_{n}+\epsilon_{n} &= \hat{g}_n - d_{n}, & \forall n \in\mathcal{N} 
\label{eq:P1b} \\
\sum_{n=1}^N q_{n} &= 0 & 
\label{eq:P1c} \\
\hat{g}_n &\leq \sum_{m \in \mathcal{M}_{n}}\left(g_{m}^{0}+r_{m}^{\uparrow}\right), & \forall n \in\mathcal{N} 
\label{eq:P1d1} \\
\hat{g}_n &\geq \sum_{m \in \mathcal{M}_{n}}\left(g_{m}^{0}-r_{m}^{\downarrow}\right), & \forall n \in\mathcal{N} 
\label{eq:P1d2} \\
\sum_{m \in \mathcal{M}_{n}}\left(g_{m}^{0}-r_{m }^{\downarrow}\right) & \geq \sum_{m \in \mathcal{M}_{n}} g_{m}^{\min } u_{m}, & \forall n \in\mathcal{N} \label{eq:P1e} \\
\sum_{m \in \mathcal{M}_{n}}\left(g_{m}^{0}+r_{m}^{\uparrow}\right) & \leq \sum_{m \in \mathcal{M}_{n}} g_{m}^{\max } u_{m}, & \forall n \in\mathcal{N} \label{eq:P1f} \\
-f_{l}^{\max } & \leq \sum_{n=1}^N h_{l n} q_{n} \leq f_{l}^{\max }, & \forall l \in\mathcal{L} \label{eq:P1g} \\
d_{n}^{\min } & \leq d_{n} \leq d_{n}^{\max }, & \forall n \in\mathcal{N} \label{eq:P1j} 
\end{align}

The objective function \eqref{eq:P1a} minimizes the ``cost'' of total power imbalance, which consists of renewable curtailment and load shedding, as captured by the slack variables $\epsilon_{n}$ and weighed by an imbalance price $\gamma$. Constraint \eqref{eq:P1b} determines the net power injection at each bus, while constraint \eqref{eq:P1c} guarantees power balance across the system. The dispatchable generation capacity limits, and up and down reserves limits are defined by constraints \eqref{eq:P1d1}--\eqref{eq:P1f}, respectively. Constraint \eqref{eq:P1g} enforces flow limits on the transmission lines. 
Finally, residual demand vector limits are enforced using minimum and maximum limits, as denoted by $d_{n}^{\min}$ and $d_{n}^{\max}$, respectively. We will refer to this model as the benchmark approach (BA) for evaluating residual demand admissibility.

\subsection{Limitations of the Benchmark Residual Demand Admissibility Assessment Approach}

The BA calculates allowable generation reserve deployment variables bus-per-bus $\hat{g}_n = \sum_{m \in \mathcal{M}_n} g_m$, net power injections $q_n$, allowable residual demands $d_n$ and their possible curtailment $\epsilon_n, \; \forall n \in \mathcal{N}$ based on the specified range of residual demands set in \eqref{eq:P1j}.

The first drawback of the BA relates to the modeling of residual demand range limits \eqref{eq:P1j}, which is a simple box in $N$-dimensional space. This constraint may lead to sub-optimal or even incorrect admissibility assessments because it ignores the spatial correlation which may exist among residual demands at different locations of a power system. Instead, we argue in favor of the use of a PUS, $d_{n} \in P(\mathcal{S})$ as described in Section~II which captures those correlations.

Another shortcoming of the BA occurs when the $N$-dimensional box of residual demands \eqref{eq:P1j} is small enough that none of the inequality constraints \eqref{eq:P1d1}--\eqref{eq:P1g} are active. In this case, the objective function of the BA equals zero, which means that the power system has sufficient reserve and transmission capacity to handle all allowable residual demand values. In these cases, one may argue that too much flexibility---in the form of reserves---was scheduled, which may be grossly uneconomical. On the other hand, there may also be cases where the margin offered by scheduled flexibility resources and available transmission capacity with respect to the range of residual demands is in fact very small. These are risky situations, where simple errors in characterizing the residual demand range, may lead to infeasibility, albeit the original flexibility assessment having come back with a positive outcome.

Clearly, the BA fails to assess the possible risks associated with poorly-assessed uncertainty that might lead to infeasibility. In this paper, we will propose a framework that has the ability to assess \emph{both} excessive and inadequate system flexibility by assessing how far operating points are from their feasibility region boundaries. We recall that the aim of flexibility assessment methods such as ours and \cite {W.Wei2, C.Wang3,C.Wang4,C.Wang2} is to estimate possible power imbalances; they do not calculate optimal reserve deployment such as dispatch problems \cite {Y.Huo1, A.Velloso}.


\subsection{Feasibility Region Projection onto the Residual Demand Space}

The first step in recasting the BA, is to determine its corresponding \emph{loadability set}. The loadability set of a power system is a set of residual demand realizations that can be supplied by generation while respecting all transmission and reserve capacity limits \cite{A.Kalantari}. A loadability set in the generation-demand space, $\Xi_{gd}(\zeta)$, is characterized as a function of $\zeta = [\zeta_1 \; \cdots \; \zeta_M]^\top$, a specific unit commitment solution 
\begin{equation}
 \Xi_{gd}(\zeta)=\left\{(g,d) \mid \eqref{eq:P1d1}-\eqref{eq:P1j} \right\}
\end{equation}
On the other hand, the loadability set, $\Xi_{d}(\zeta)$, is the projection of $\Xi_{gd}(\zeta)$ onto demand space only. It is defined as
\begin{equation}
 \Xi_{d}(\zeta)=\left\{d \mid \exists(g, d) \in \Xi_{gd}(\zeta) \right\}
\end{equation}

The loadability set $\Xi_{gd}(\zeta)$ is a relaxation of the BA original feasibility region. It enforces all its inequalities, but ignores its equality constraints. Alternatively, $\Xi_{d}(\zeta)$ describes which residual demand vectors can be supported by a power system considering a specific unit commitment instance $\zeta$. As seen in \cite{A.A.Jahromi1} and \cite{A.A.Jahromi2}, the projection process whereby we pass from $\Xi_{gd}(\zeta)$ to $\Xi_{d}(\zeta)$ involves the generation of large numbers of constraints, many of which are redundant.

To counter the explosion in the number of redundant constraints associated with the mapping of $\Xi_{gd}(\zeta)$ onto $\Xi_{d}(\zeta)$, we adopt the iterative approach for determining the minimal representation of the loadability proposed by \cite{A.A.Jahromi2}. This approach combines the \emph{umbrella constraint discovery} (UCD) algorithm and the Fourier–Motzkin elimination (FME) method as follows. First, we identify network constraints which are \emph{umbrella}---i.e., which are non-redundant in the BA and effectively contribute in shaping its feasibility region---, and we rule out the corresponding redundant network constraints. Second, we construct the loadability set in the generation-demand space considering the dispatchable generators capabilities and umbrella network constraints only. The loadability set constraints are mapped from the generation-demand space onto the demand space by removing the generator variables one by one using FME \cite{A.A.Jahromi1}. After each elimination, redundant constraints are identified and removed to keep the number of constraints as small as possible. We terminate the procedure once all dispatchable generator variables have been eliminated.

\section{Inverse Optimization Applied to Explicit Flexibility Characterization}\label{METHODOLOGY} 

In this section, inspired by recent advances in data-driven inverse optimization techniques \cite{T.C.Y.Chan2, K.Ghobadi, Aaron}, we present how inverse optimization can be used to address the flexibility assessment problem. Inverse optimization describes the ``reverse'' process of the conventional mathematical optimization. An inverse optimization problem takes decisions or observations as input and determines the objective function and/or the constraints that render these observations approximately or exactly optimal.

\subsection{Generalized Inverse Linear Optimization Problems}
Here, we introduce a general inverse optimization (GIO) model \cite{Aaron} for a linear optimization problem without any explicit assumptions regarding its feasibility as in \cite{T.C.Y.Chan2}. First, we start from a linear optimization problem which is called the \emph{forward problem} (FO). Let ${c,x} \in \mathbb{R}^{N}$ denote cost and decision vectors, respectively. Let ${A} \in \mathbb{R}^{J \times N}$, and ${b} \in \mathbb{R}^{J}$ denote the constraint and the right hand side vector, respectively.\footnote{We are using dimensions $J$ and $N$ in a generic sense here; at this stage, there is no explicit connection with prior developments presented earlier in the paper.}

\begin{equation}
\label{objective_Minimum}\
\mathrm{FO}\left(c\right):
\min c^{\top}{x}
\end{equation}
Subject to:
\begin{equation}
A x \geq b 
\end{equation}

Assuming that the forward problem does not have redundant constraints, given a decision $\hat{{x}} \in \mathbb{R}^{N}$ and for an $r$ norm $(r \geq 1)$, the single observation \emph{generalized inverse linear optimization problem} (GIO) is
\begin{equation}
\mathrm{GIO}_r \left(\hat{{x}}\right): \min_{c, y, s} \|s\|_r \label{eq:P3a}
\end{equation} 
Subject to:
\begin{equation}
 \quad A^{\top} y= c, \quad y \geq 0 \label{eq:P3b}
\end{equation}
\begin{equation}
c^{\top} \hat{x} = b^{\top} y+c^{\top} s \label{eq:P3c}
\end{equation}
\begin{equation}
A\left(\hat{x}-s\right) \geq b \label{eq:P3d}
\end{equation}
\begin{equation}
\|c\|_{1}=1 \label{eq:P3e}
\end{equation}
The objective function \eqref{eq:P3a} minimizes the error (perturbation) vector $s \in \mathbb{R}^{N}$ using an arbitrary $r$-norm, which provides a natural measure of error in the space of decision variables. In GIO, $y \in \mathbb{R}^{J}$ represents the vector of dual variable of the forward problem's constraints. Constraints in \eqref{eq:P3b} enforces the dual feasibility of the forward problem. Constraint \eqref{eq:P3c} connects the cost vector of the forward problem and its dual variables' vector, including the perturbation vector whose value may have to be non-zero to satisfy strong duality of the forward problem for the given $\hat{x}$. Constraint \eqref{eq:P3d} enforces primal feasibility of the perturbed decisions $\hat{x}-s$. Finally, constraint \eqref{eq:P3e} normalizes of the cost vector to prevent it from collapsing to the trivial solution, which is zero.


\subsection{GIO for Flexibility Assessment}

Consider a power system's feasible space whose minimum realization loadability set has been determined ($\mathcal{D}(\zeta)$) based on \eqref{eq:P1d1}--\eqref{eq:P1j}, or, more accurately, by imposing $d \in P(\mathcal{S})$ in lieu of \eqref{eq:P1j}. We demonstrate next how GIO can assess the ``distance'' of any $d$ to the boundaries of $\mathcal{D}(\zeta)$. We posit that knowledge of how far (or near) an operating point is to the boundary of its feasibility region is highly valuable to power system operators and planners alike. With such knowledge at hand, one could seek to extend that distance for increased robustness or try to reduce it, when excessively large, to aim for more economical system operations.

We now consider some instance $d^{0}$ (as $\hat{x}$ in the general case), as well as vectors $\alpha'_j$ and $\beta'_j$, $j = 1, \ldots, J'$---exogenously determined by calculating $\mathcal{D}(\zeta)$---vectors to form the input data set to infer GIO variables. Aiming for a more compact notation, we stack the $J'$ inequalities of $\mathcal{D}(\zeta)$ into an $\mathbb{R}^{J'\times N}$ matrix $A$ (whose rows are $(\alpha'_j)^\top, \; \forall j \in \mathcal{J}'$) and a $\mathbb{R}^{J'\times 1}$ vector $b$ (whose rows are $\beta'_j, \; \forall j \in \mathcal{J}'$); from these we now have $\mathcal{D}(\zeta) = \{ d \in \mathbb{R}^N \mid A d \geq b \}$.

Therefore, we define the data-driven inverse optimization problem (DDIO) for flexibility assessment
\begin{equation} \label{eq:P4a}
\mathrm{DDIO}_r \left( d^{0},\zeta \right): \min_{c,y,s} \|{s}\|_r
\end{equation}
Subject to:
\begin{equation}\label{eq:P4b}
A^{\top} {y}=c , \quad {y} \geq 0 
\end{equation}
\begin{equation}\label{eq:P4c}
c^{\top}\left(d^{0}-{s}\right)=b^{\top} {y} 
\end{equation}
\begin{equation}\label{eq:P4d}
A\left(d^{0}-{s}\right) \geq b
\end{equation}
\begin{equation}\label{eq:P4e}
\|c\|_{1}=1
\end{equation}

When the objective function of DDIO is equal to zero, this means $d^{0}$ is optimal. In the case where $\mathcal{D}(\zeta)$ is a polytope and if $\| s \|_r = 0$, it effectively means that $d^0$ is one of the vertices of $\mathcal{D}(\zeta)$. This is so because the optimal solutions of linear programs on bounded convex polyhedrons (polytopes)---like $\mathcal{D}(\zeta)$---are strictly located on their vertices. In fact, it is possible to construct a $c$ that puts $d^0$ on a vertex of $\mathcal{D}(\zeta)$, that is where bounded linear programming problem solutions lie.\footnote{This interpretation is specific to linear costs $c^\top x$ only. If the assumed cost function was nonlinear (e.g., quadratic), the solution could be at any point on the boundary of $\mathcal{D}(\zeta)$.} Otherwise, $d^0$ is interior or exterior to $\mathcal{D}(\zeta)$. Therefore, the vector $s \neq 0$, and the norm of $s$ can be interpreted as a measure of distance to the closest boundary of $\mathcal{D}(\zeta)$.

We note that DDIO \eqref{eq:P4a}--\eqref{eq:P4e} is a nonlinear, non-convex optimization problem because of the of the $c^\top s$ term in the strong duality constraint \eqref{eq:P4c}. In the next subsection, we will propose a computationally efficient method to solve DDIO in reasonable time using off-the-shelf optimization software.

\begin{figure}
\centering
\includegraphics[width=6.5cm]{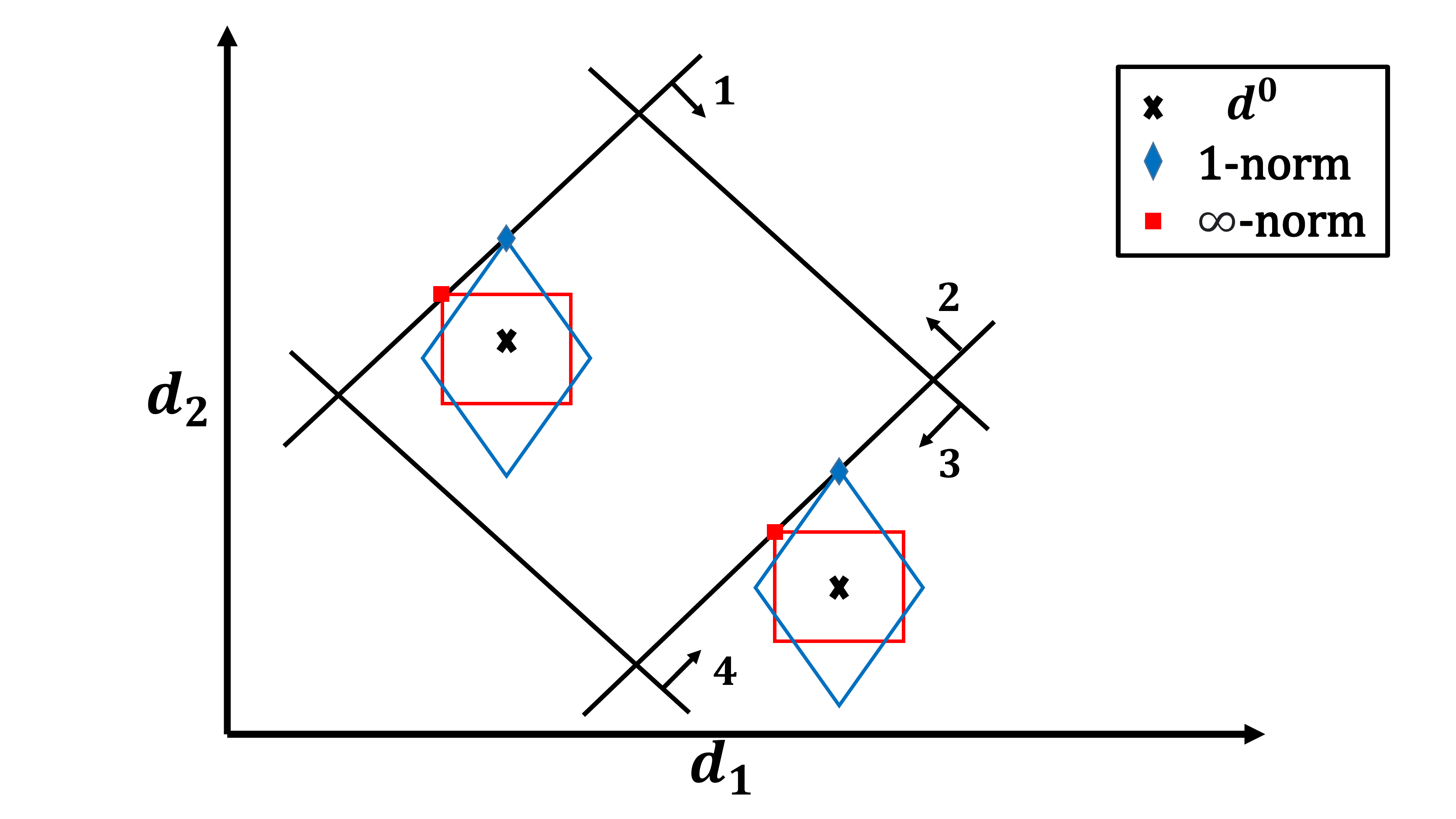}
\caption{Geometrical interpretation of inverse optimization solution methodology.}\label{fig:02}
\end{figure}

To overcome the non-convexity of \eqref{eq:P4c}, we leverage the feasibility projection problem structure described in \cite{K.Ghobadi, Aaron} to solve the inverse problem. The problem can be decomposed into $j = 1,\ldots,J'$ linear sub-problems
\begin{equation} \label{eq:P5a}
\min_{s_j} \left\|s_j\right\|_{r}
\end{equation}
Subject to:
\begin{equation} \label{eq:P5b}
A(d^{0}-s_{j}) \geq b
\end{equation} 
\begin{equation}\label{eq:P5c}
\left(\alpha'_{j}\right)^{\top} (d^{0}-s_{j}) = \beta'_{j}
\end{equation}

and then picking out the smallest $\|s_j\|_r$ over all $j \in \mathcal{J}'$. The problem \eqref{eq:P5a}--\eqref{eq:P5c} finds the shortest distance between the residual demand forecast $d^0$ and each of the constraints $j \in \mathcal{J}'$. The vector $s_j$ that is found has to be consistent with all constraints \eqref{eq:P5b}, and it has to bring $d^0$ in contact with constraint $j$ as per \eqref{eq:P5c}. Picking out the constraint or constraints $j^* \in \mathcal{J}^* \subseteq \mathcal{J}'$ with the shortest distance finds the closest point to the boundary of $\mathcal{D}(\zeta)$ with respect to $d^0$. A geometrical illustration of \eqref{eq:P5a}--\eqref{eq:P5c} is shown in Fig.~\ref{fig:02}. Here $d^0$ can be in the interior of $\mathcal{D}(\zeta)$, or its exterior, and, based on the specific choice of a vector norm, we can determine how far $d^0$ is from its boundary. 

\subsection{Flexibility Assessment Through DDIO}
 
As just seen, DDIO is capable of determining the shortest distance between a residual demand vector $d^0$ and the boundaries of a power system's loadibility set $\mathcal{D}(\zeta)$. Moreover, it is able to identify which constraint or constraints of  $\mathcal{D}(\zeta)$ are either violated (if $d^0$ is exterior) or close to be violated (if $d^0$ is interior). Therefore, loosely speaking, a residual demand forecast is ``flexibility adequate'' if $d^0$ is in the interior of $\mathcal{D}(\zeta)$ and has a large distance to its boundaries. 

Based on these observations, we thus propose flexibility metrics based DDIO results:
\begin{itemize}
\item Analytical flexibility metric $\rho$: Inspired by\cite{T.C.Y.Chan2}, we argue that
\begin{equation}
 \rho_{r}\left(d^{0}, \zeta \right)
= 1-\frac{\left\|s_{j^*}\right\|_{r}}{(1 /J') \sum_{j=1}^{J'}\left\| s_j\right\|_{r}}
\end{equation}
is a useful indicator of the relative security of the residual demand $d^0$. It weighs the shortest distance between $d^0$ to the boundaries of $\mathcal{D}(\zeta)$ relative to the average of boundary distances. The indicator ${\rho}_r \to 1$ if $d^0$ is close to one or more boundaries (i.e., system flexibility is in tight supply and/or may not be transmittable to all portions of the network), or it would tend to zero if it is somewhat in the center of $\mathcal{D}(\zeta)$ (i.e., system flexibility is plenty and can be adequately transmitted across the network). We distinguish one weakness of $\rho_{r}$, however. As formulated here, it is not possible to distinguish whether or not $s_j$ is used to bring $d^0$ to a boundary of $\mathcal{D}(\zeta)$ from the inside or the outside. In fact, we would argue that $\rho_{r}$ is useful when $d^0 \in \mathcal{D}(\zeta)$. For cases where $d^0 \notin \mathcal{D}(\zeta)$, we propose the next metric.


\item Residual demand curtailed (RDC): In case there exists one or more $j \in \mathcal{J}'$ such that $(\alpha'_j)^\top d^0 < \beta'_j$ (i.e., $d^0 \notin \mathcal{D}(\zeta)$), residual demand may need to be curtailed as indicated by the components of $s_j$. Defining the subset of violated boundaries of $\mathcal{D}(\zeta)$, $\tilde{\mathcal{J}} = \{ j \in \mathcal{J}' \mid (\alpha'_j)^\top d^0 < \beta'_j \}$, we have
\begin{equation}\label{RDC}
\mbox{RDC} = \sum_{j \in \tilde{\mathcal{J}}} \sum_{n = 1}^N s_{jn}
\end{equation}
where we note that if RDC~$> 0$ demand has to be shed, and if RDC~$< 0$ non-dispatchable generation has to be curtailed.
\end{itemize}

We will see in the next section how the choice of norm can influence values of $\rho_{r}$ and RDC. Moreover, the advised reader may be tempted to suggest computing the volume of the loadability set as a way to assess flexibility adequacy. As exposed in \cite{J.Lawrence}, this would be prohibitive from a computational point of view. The computation effort required to obtain the volume of a polytope increases very rapidly with its number of vertices. In addition, we recall from \cite{A.A.Jahromi1} that the number of vertices of loadability sets tends to grow very rapidly with successive applications of  Fourier-Motzkin elimination. Therefore, it would be ill-advised to pursue this objective here. In our case studies presented next, we will nonetheless compute the volume of $\mathcal{D}(\zeta)$ for comparison purposes; however, we shall use an approach based on Monte Carlo simulations rather than compute volumes analytically \cite{F.Golestaneh}. 

\section{Case Studies}\label{Case}
In this section, we illustrate how the use of forward-looking residual demand PUS and its coupling with DDIO can work to assess post-unit commitment flexibility adequacy.

\subsection{Procedure for Constructing Correlated Residual Demand Forecast Error Time Series}

We generate $N$ synthetic spatially-correlated residual demand time series of length $T$, which are then consigned to matrix $W$. They consist of historic residual demand forecasts $\mu \in \mathbb{R}^{T \times N}$, which correspond to the nominal demand values from the data sets in \cite{Data1}. These are superimposed with zero-mean normally-distributed forecast errors with spatial correlation given by the covariance matrix $\Sigma$. Here, we take the approach outlined in \cite{Pena}, where errors are assumed to be proportional to forecasts and whose variance and correlation are adjusted with an uncertainty level parameter $\eta \in [0,1]$. This way, the diagonal elements of $\Sigma$ are $\sigma^2_{nn} = (\eta \mu_n)^2$ for all $n \in \mathcal{N}$, and, as proposed in \cite{Pinar}, its off-diagonal elements are given by $\sigma^2_{nn'} = \eta^2 \alpha \mu_n \mu_{n'}$, for all $n \neq n' \in \mathcal{N}$, and where $\alpha \in [-1,1]$ is an adjustable correlation coefficient.

\subsection{Three-Bus Test System}
Here we use the three-bus test system from \cite{A.A.Jahromi1} with a slight modification by adding two wind farms at buses 2 and 3, and neglecting the two generators ramping constraints from \cite{A.A.Jahromi2}. We consider two scenarios whose past forecasted residual demands were 320 and 50~MW (Scenario 1) and 240 and 40~MW (Scenario 2) at buses 2 and 3, respectively. We generated $T=4000$ time instances for each scenario as seen in Fig. 4 (a) for Scenario 1. In Scenario 1, residual demands have an uncertainty level $\eta$ of 0.067 per unit and a correlation coefficient $\alpha$ of 0.8. For Scenario 2, $\eta = 0.1$ per unit with $\alpha = 0.7$. We assume that the penalty for load shedding and renewable generation spillage is $\gamma =\$1000$ per MWh.

We define three intact loadability sets which do not assume residual demand forecast errors and its uncertainty; that is, they exclude \eqref{eq:Pus1} and \eqref{eq:P1j} and only require that residual demands be strictly positive. These sets are called $\mathcal{D}(\zeta_{1}) = \mathcal{D}(u_1=1,u_2=0)$, $\mathcal{D}(\zeta_{2}) = \mathcal{D}(u_1=0,u_2=1)$ and $\mathcal{D}(\zeta_{3}) = \mathcal{D}(u_1=1,u_1=1)$, and they are described by six, five, and six constraints, respectively. Moreover, for these we constructed loadability sets based on both PUS and box uncertainty sets, and we refer to them as PUS$\mathcal{D}$ and Box$\mathcal{D}$, respectively.

\begin{figure}
\centering
\begin{subfigure}{0.23\textwidth}
    \includegraphics[width=\textwidth]{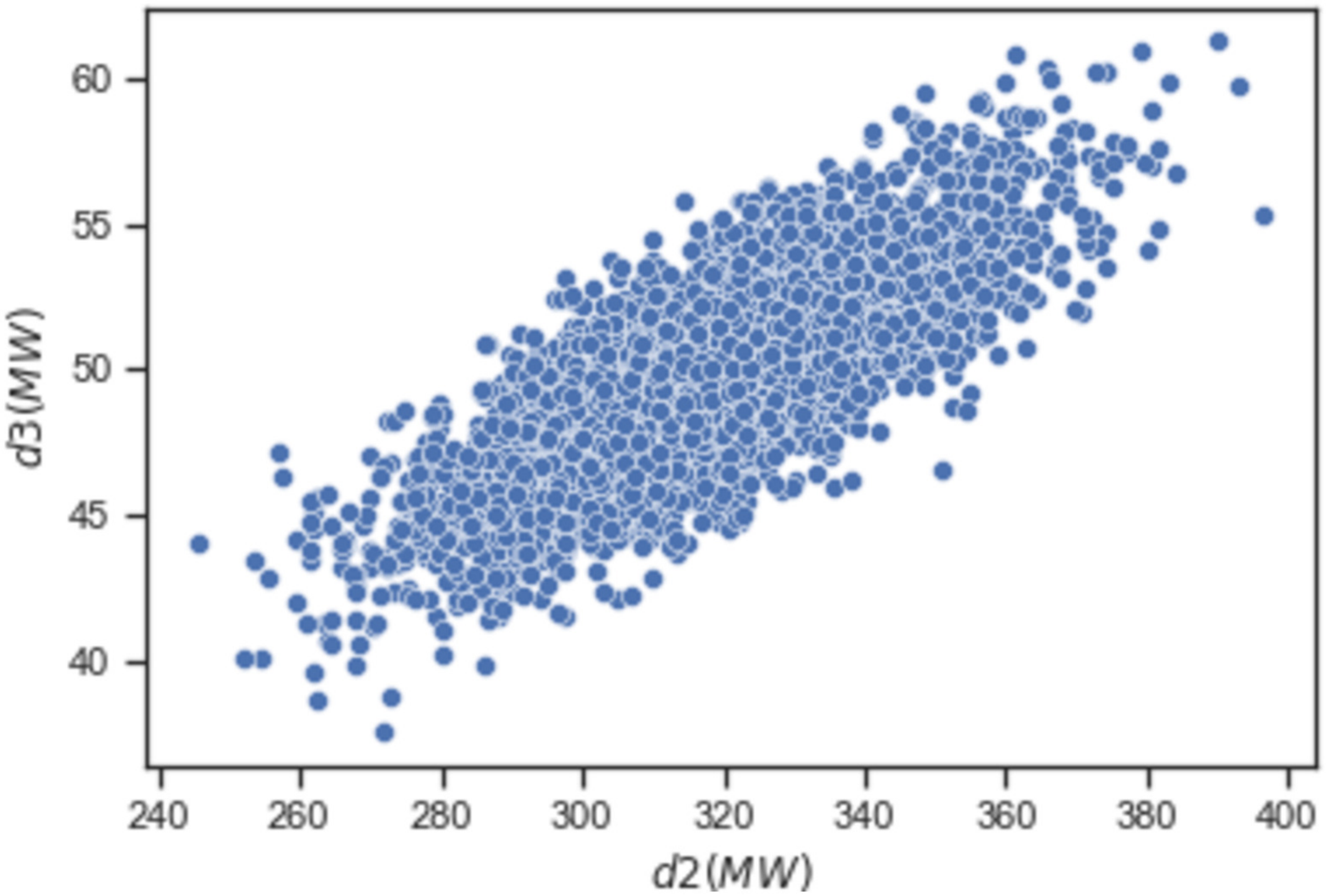}
    \caption{}
    \label{fig:04a}
\end{subfigure}
\hspace{0.5em}
\begin{subfigure}{0.23\textwidth}
    \includegraphics[width=\textwidth]{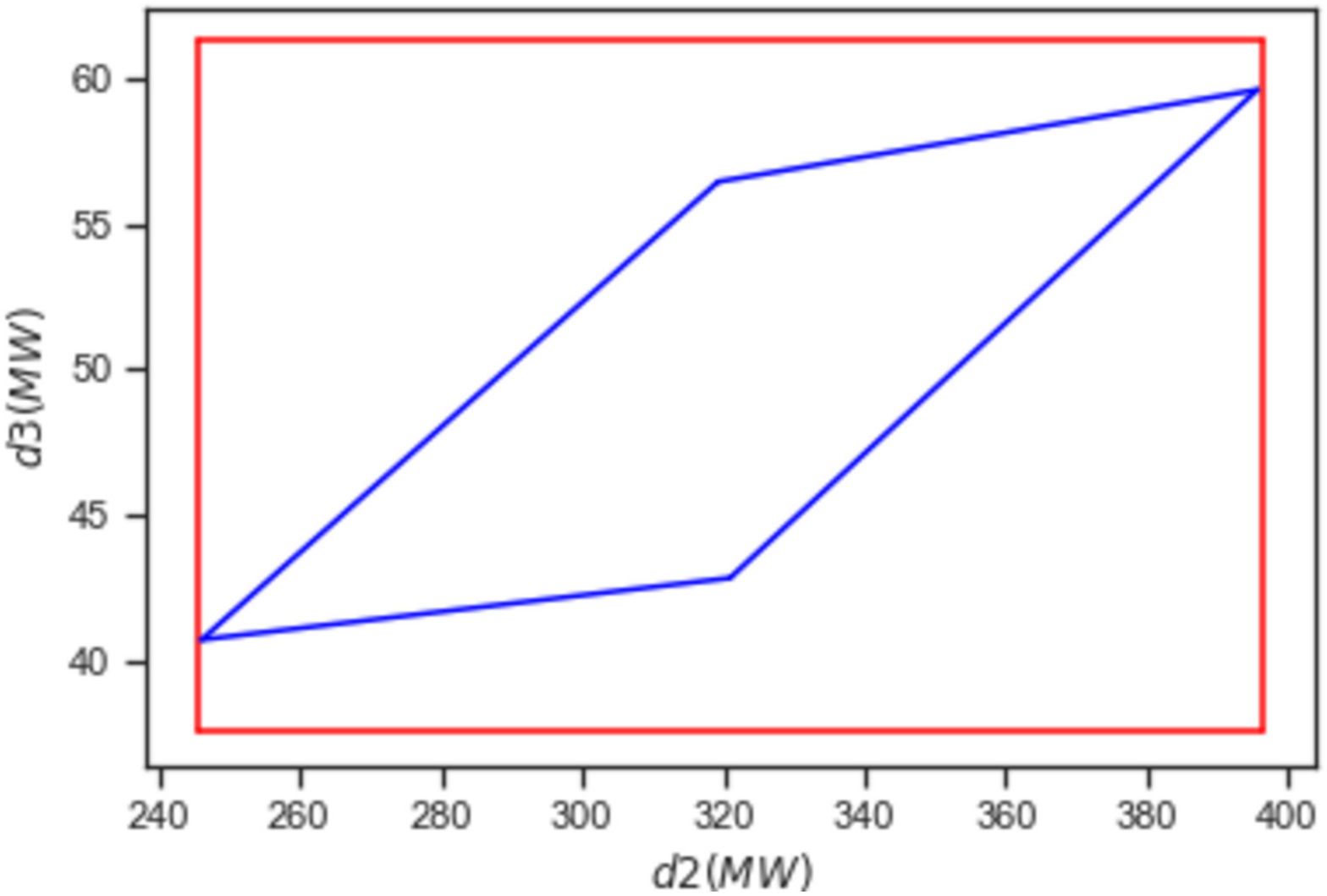}
    \caption{}
    \label{fig:04b}
\end{subfigure}
\hfill
\caption{(a) Scatter plot of residual demands for Scenario 1; (b) Uncertainty sets for Scenario 1.}
\label{fig:04}
\end{figure}

\begin{table}
\captionsetup{font=footnotesize}
\caption{\sc{Uncertainty and Loadability Set Volumes}} 
\centering 
\footnotesize{
\begin{tabular}{c c c c c} 
\hline
\multicolumn{5}{c}{Volume (MW$^{2}$)}\\ [0.1ex] 
\hline
& \multicolumn{2}{c}{Uncertainty set} & \multicolumn{2}{c}{Loadability set} \\ [1.0ex] 
\hline
Scenario & PUS &  Box  & PUS$\mathcal{D}(\zeta_{3})$ & Box$\mathcal{D}(\zeta_{3})$ \\ [0.1ex]
\hline
1 & 1163 &  3950   &   936   & 3520 \\ [0.1ex]
2 & 1929 &  4740   &   1810  & 4530 \\ [0.1ex]
\hline 
\end{tabular}}
\label{table1}
\end{table}

\begin{table}
\captionsetup{font=footnotesize}
\caption{\sc{Performance of Flexibility Assessment Strategies}} 
\centering 
\footnotesize{
\begin{tabular}{c c c c}
\hline
Scenario & Method & RDC  \\\hline
\multirow{3}{*}{1} & $\operatorname{DDIO}^{1}$ & 98.527 (0.0\%) \\
& $\operatorname{DDIO}^{\infty}$ & 147.48 (+50\%) \\
& BA & 98.527 \\\hline 
\multirow{3}{*}{2} & $\operatorname{DDIO}^{1}$ & 39.065 (0.0\%) \\2
& $\operatorname{DDIO}^{\infty}$ & 39.065 (0.0\%) \\
& BA & 39.065 \\
\hline 
\end{tabular}}
\label{table2} 
\end{table}

\maketitle
\begin{figure*}%
\centering
\begin{subfigure}{.6\columnwidth}
\includegraphics[scale=0.16]{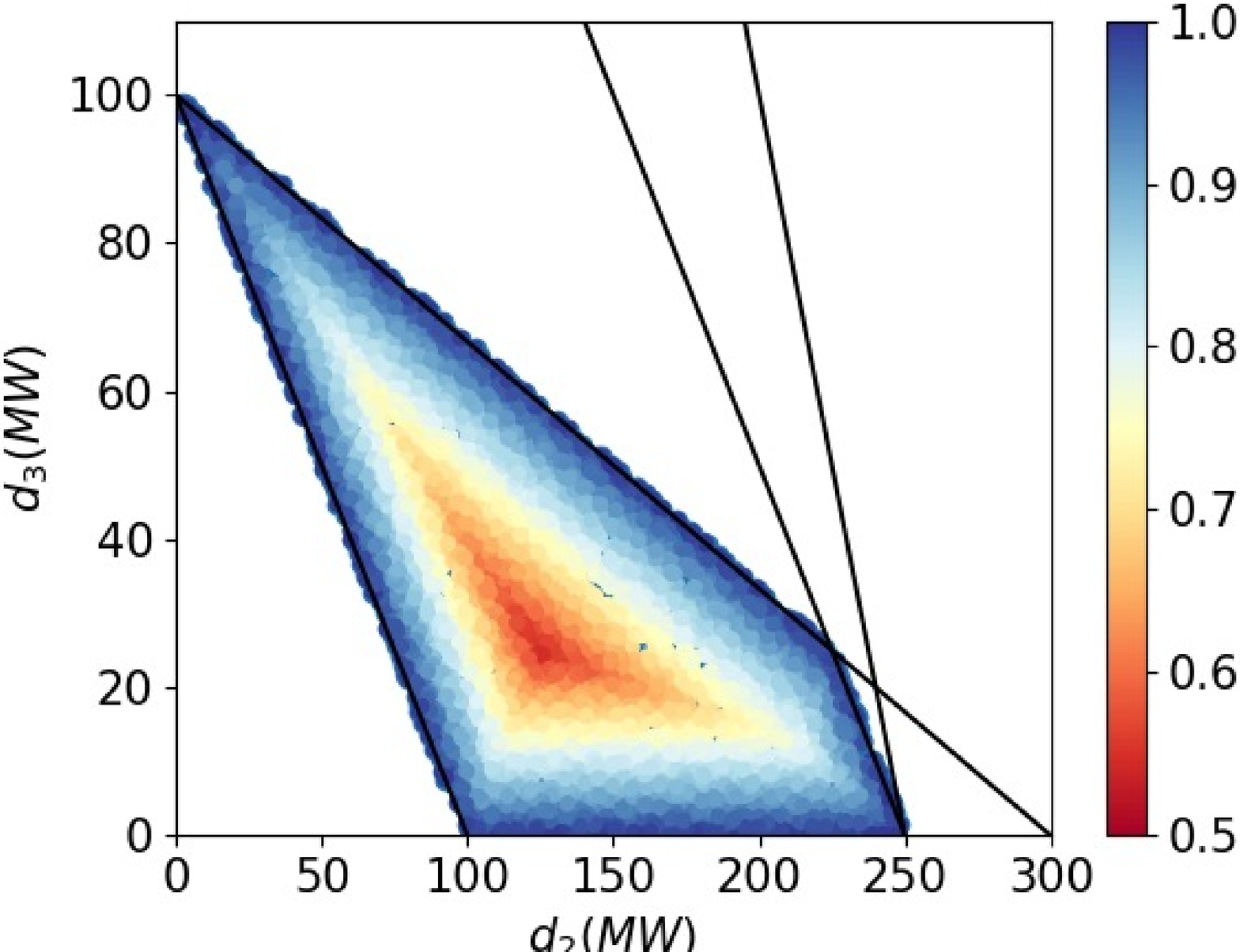}%
\caption{$\mathcal{D}(\zeta_{1})$}%
\label{subfig05a}
\end{subfigure}\hfill%
\begin{subfigure}{.6\columnwidth}
\includegraphics[scale=0.16]{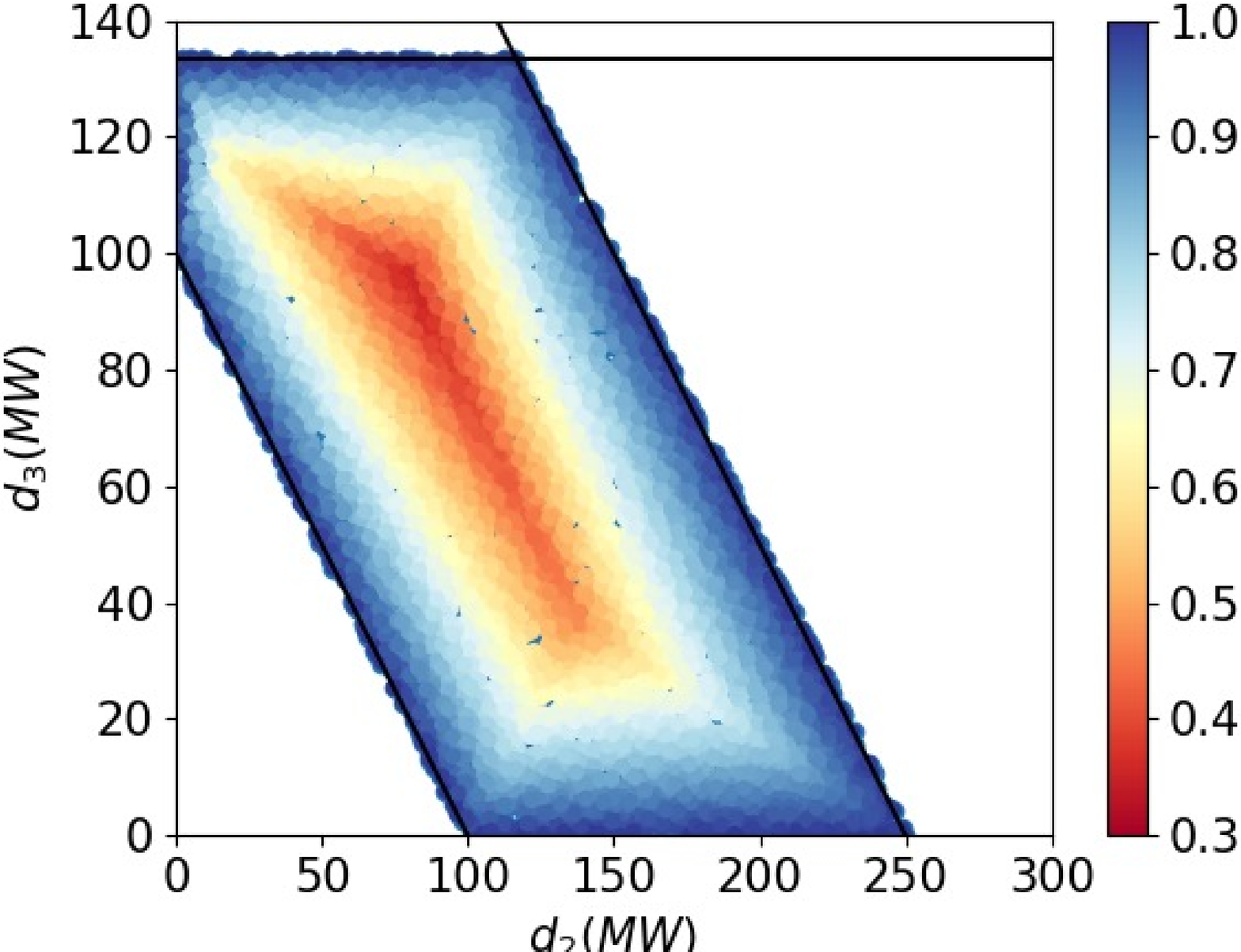}%
\caption{$\mathcal{D}(\zeta_{2})$}%
\label{subfig05b}
\end{subfigure}\hfill%
\begin{subfigure}{.6\columnwidth}
\includegraphics[scale=0.16]{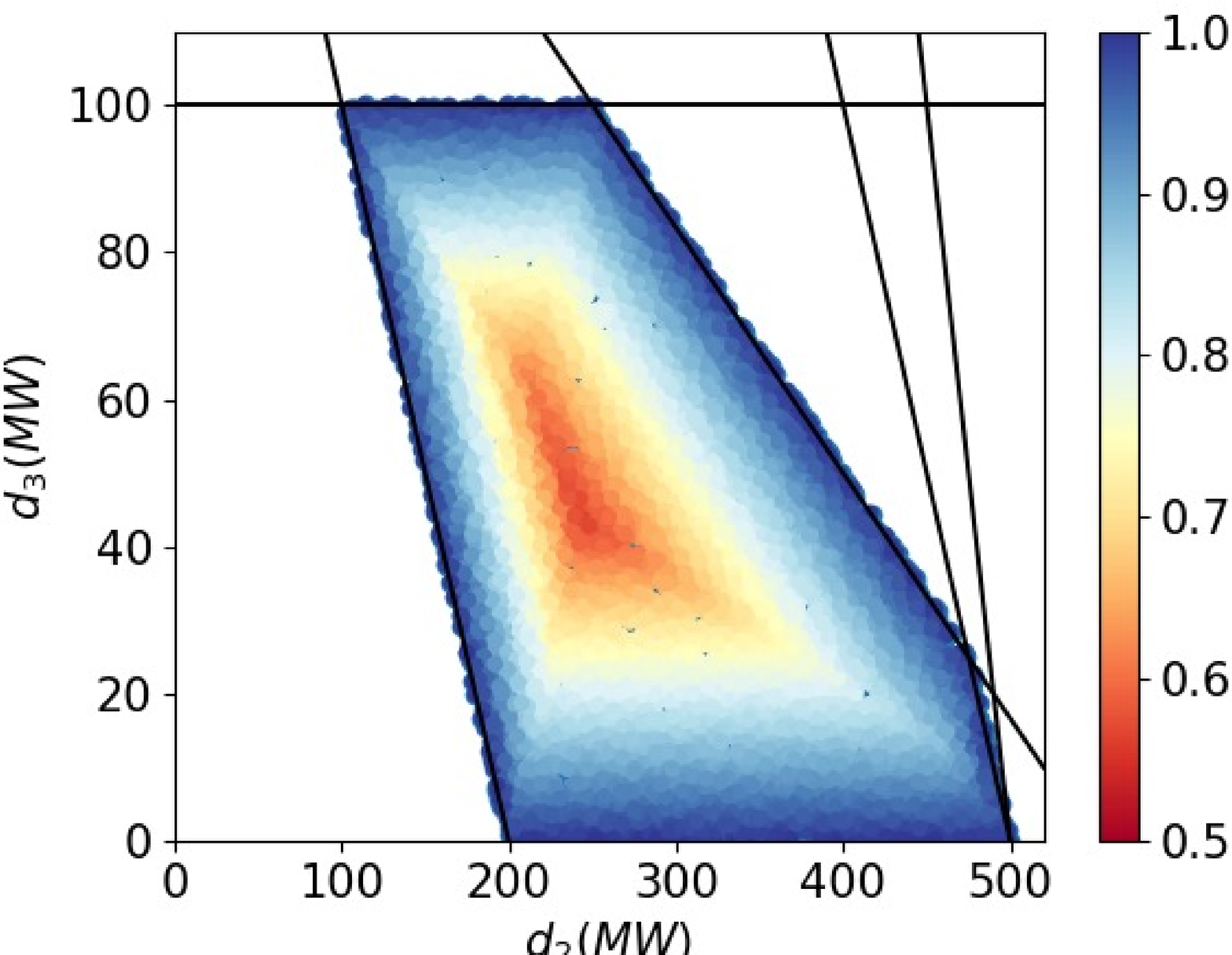}%
\caption{$\mathcal{D}(\zeta_{3})$}%
\label{subfig05c}
\end{subfigure}%
\caption{Flexibility metric ${\rho}_{\infty}$ distribution across the intact loadability sets.}
\label{fig:05}
\end{figure*}

\subsubsection{Uncertainty and Loadability Sets' Volume}

In Fig.~\ref{fig:04} (b), one can see the PUS for the data from Scenario 1 (blue rhombus); in addition, we provide the box set corresponding to \eqref{eq:P1j} (red rectangle). By inspection, we see that the PUS can effectively capture the spatial correlation between the residual demands, while being less conservative than the boxed uncertainty set. In Table.~\ref{table1}, we provide the volume (areas here) for the two scenarios considering PUS and box-shaped uncertainty modeling. The conservative box uncertainty set is 3.4 and 2.45 larger than PUS for Scenarios 1 and 2, respectively.

Also, in Table~\ref{table1} we characterize the loadability sets considering two generators, the  network constraints and data-driven uncertainty sets. For Scenario 1, the loadability sets PUS$\mathcal{D}(\zeta_{3})$ and Box$\mathcal{D}(\zeta_{3})$ are characterized by their uncertainty boundaries plus one constraint representing the umbrella line flow constraint after eliminating generators variables using FME and removing any redundancy using UCD. Similarly for Scenario 2, PUS$\mathcal{D}(\zeta_{3})$ and Box$\mathcal{D}(\zeta_{3})$ are characterised by their uncertainty sets and one constraint representing the boundary defined by the generator minimum capacity limits. Clearly, system constraints associated to respective loadability sets shave larger volumes from the box set in comparison to PUS in both scenarios. In fact, Box$\mathcal{D}(\zeta_{3})$ volume's was shaved off by 1.90 and 1.76 times more than PUS$\mathcal{D}(\zeta_{3})$ for Scenarios 1 and 2, respectively.

\subsubsection{Flexibility Assessment Strategies' Performance}

Table~\ref{table2} displays the performance of three proposed methods for both scenarios and the loadability set $\mathcal{D}(\zeta_{3})$. The percentage changes in RDC are taken with respect to the BA. Both scenarios emphasize that $\operatorname{DDIO}^{1}$ (DDIO solved using the $r=1$ norm) is able to assess residual demand curtailed in the same way as the BA. In the case of Scenario~1, $\operatorname{DDIO}^{\infty}$ shows an increase in RDC of 50\% compared to $\operatorname{DDIO}^{1}$ and the BA. The reason for this is that when using the $r={\infty}$ norm for calculating the shortest distance to reach the loadability set boundary from $d^0$, DDIO is using all dimensions to get to the boundary. On the other hand, when the $r={1}$ norm is used, DDIO is using only one dimension, while neglecting the others. However, in the calculation of RDC in \eqref{RDC}, all $N$ components of $s_j$ are used, which the the case in the calculation of RDC and in the solution of $\operatorname{DDIO}^{1}$. 




\subsubsection{Impacts of Generation Schedules on Flexibility Metrics}
Here, the properties of the flexibility metric $\rho_{r}$ are discussed under the three different loadability sets, $\mathcal{D}(\zeta_i); i = 1, 2, 3$. The variations of $\rho_{r}$ under different generating conditions are shown in Fig.~\ref{fig:05}. For $r = \infty$, we can see how ${\rho}_\infty$ varies for the three loadability sets as a unit-free quantity in the range $[0,1]$. By inspection, the flexibility metric varies in a structured way inside the loadability set. It shows a non decreasing trend in a radial fashion as the loading moves from the interior towards the loadability set boundaries. Hence, for any $d^0$ and generation schedule $\zeta$ it is possible to obtain a normalized ``flexibility score'' indicating its degree of safety in terms of relative distance to the boundaries of loadability sets.

\subsection{IEEE Reliability Test System}
Next, we test the proposed approach on the IEEE Reliability Test System (RTS). The data of the network, the nominal demand profile at each node, and the minimum and maximum power outputs of generators are adopted from MATPOWER \cite{Data1}. Generating units ramp rate constraints are neglected as in \cite{A.A.Jahromi2}. We generated $T = 4000$ time instances of residual demand while varying uncertainty levels to generate synthetic polyhedral and box uncertainty sets. Lastly, we assume here that $\zeta$ corresponds to the RTS unit commitment where \emph{all} its units are online.

\subsubsection{Loadability Set Construction}
Here, we consider a case where residual demand is uncertain at 17 of the 24 network nodes. Historical forecast errors are generated considering nominal loading levels with $\eta = 0.067$ and $\alpha = 0.7$. Polyhedral and box uncertainty sets are built to capture residual demand limits which are in turn used in characterizing the loadability sets of the RTS including its transmission limits.

We construct the loadability set constraints \eqref{eq:P1d1}--\eqref{eq:P1g} of the network along the polyhedral and box uncertainty sets in the generation-demand space considering all 10 generator nodes. For the polyhedral uncertainty set, we feed the UCD algorithm with corresponding residual demand limits along with the generation and network constraints as follows. First, we group uncertain residual demands into three separate polyhedral sets: $P_1(\mathcal{S}, d_{1-6}),P_2(\mathcal{S}, d_{7-10,13-14}), \text{and} \, P_3(\mathcal{S}, d_{15-16,18-20})$. Set $P_1(\mathcal{S}, d_{1-6})$ captures the uncertainty in residual demand seen at nodes 1--6, set $P_2(\mathcal{S}, d_{7-10,13-14})$ the uncertainty at nodes 7--10, 13 and 14, and finally $P_3(\mathcal{S}, d_{15-16,18-20})$ maps uncertainty at nodes 15, 16, 18--20 (for the grand total of 17 nodes). The construction of three separate PUS rather than a single one is driven by the need to limit the growth in the number of superfluous constraints when the UCD algorithm is run. Hence, instead of having $2^{17}$ constraints considered simultaneously, much fewer constraints need to be assessed at the UCD stage: $2^{6} \; (\text{for } P_1) + 2^{6} \; (\text{for } P_2) + 2^{5} \; (\text{for } P_3) = 160$ native PUS constraints. We choose to keep all principal components of each set to ensure robustness of the resulting polyhedral uncertainty sets.           


While there are ten generators in the RTS, we need to project only four generation variables at a time using FME to characterize each loadability set. This is done while the remaining generators are pushed to one of their capacity limits with the aim of maximizing the resulting loadability set volumes. This is a valid assumption since the maximum number of potentially active transmission constraints in the RTS is four, as shown in \cite{A.A.Jahromi2}. Going with the approach taken in \cite{A.A.Jahromi2}, we assume the marginal units are located at buses 1, 7, 16 and 22. The four marginal units are eliminated one by one using FME, and after each elimination, the UCD algorithm removes superfluous constraints. Table~\ref{table:03} lists the numbers of remaining umbrella constraints after each marginal unit elimination for both the PUS and box uncertainty sets. We can see that the number of constraints required by the PUS uncertainty is greater than that needed by the box set. This is due to the fact that the PUS is much more surgical in bounding the uncertainty.


  
\begin{table}
\captionsetup{font=footnotesize}
\caption{\sc{Number of umbrella constraints describing the loadability sets by uncertainty set type}} 
\centering 
\footnotesize{
\begin{tabular}{ccc} 
\hline
& \multicolumn{2}{c}{No. of constraints identified} \\ [1.0ex] 
\hline
 & PUS$\mathcal{D}(\zeta)$ & Box$\mathcal{D}(\zeta)$  \\ [1.0ex]
\hline 
In generation-demand space        & 209    & 87    \\[0.3ex]
After eliminating $\hat{g}_1$     & 200    & 77    \\[0.3ex]
After eliminating $\hat{g}_7$     & 180    & 55    \\[0.3ex]
After eliminating $\hat{g}_{16}$  & 178    & 53     \\[0.3ex]
After eliminating $\hat{g}_{22}$  & 176    & 51     \\[0.3ex]
In demand space                   & 161    & 35     \\[0.3ex]
\hline
\end{tabular}}
\label{table:03} 
\end{table}

\begin{figure}
\centering
\includegraphics[width=8.5cm]{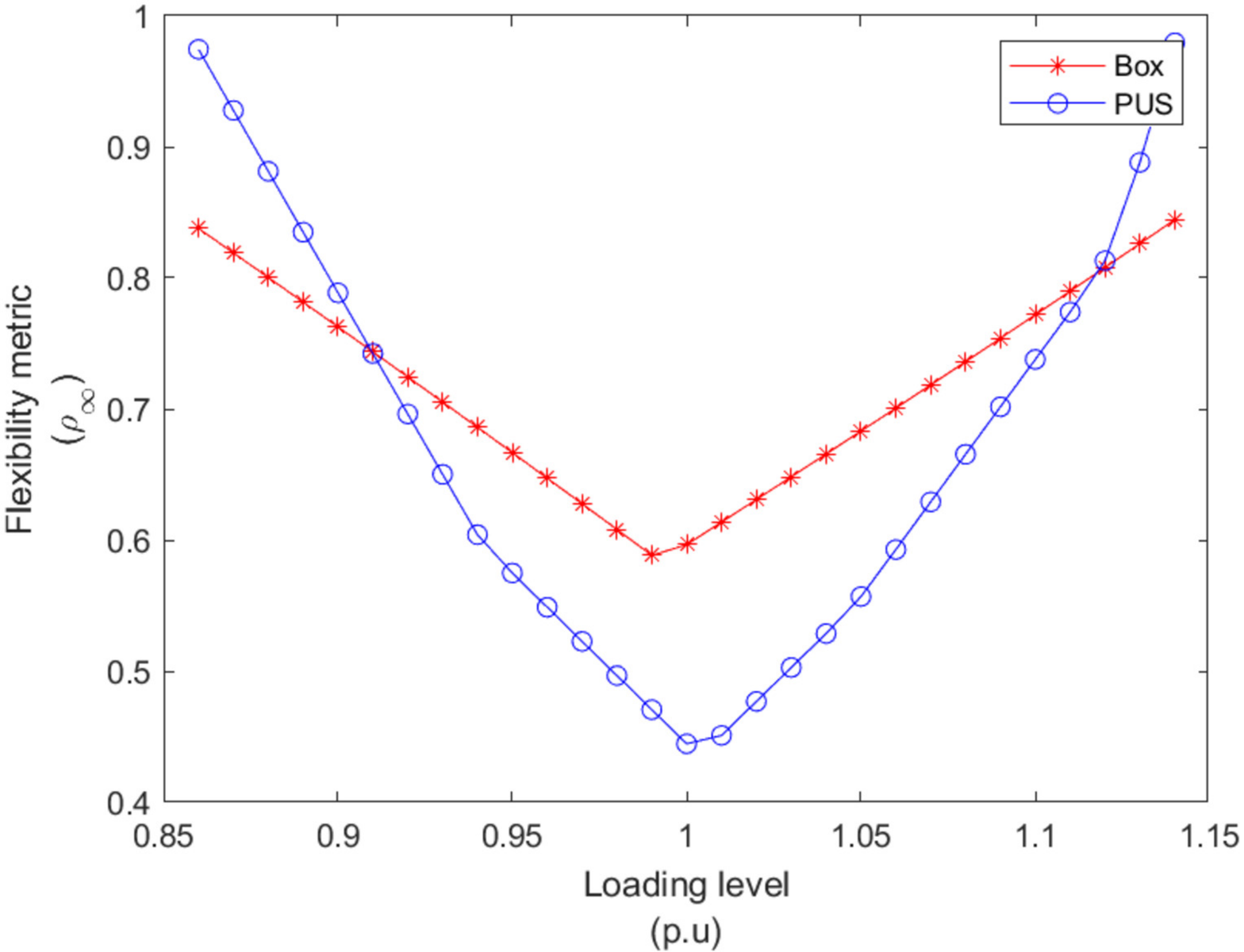}
\caption{Flexibility metric under PUS$\mathcal{D}(\zeta)$ and Box$\mathcal{D}(\zeta)$.}\label{fig:06}
\end{figure}

\subsubsection{Flexibility Assessment Strategies' Performance}
Fig.~\ref{fig:06} displays the flexibility metric $\rho_\infty$ while varying the loading  in the range of $\pm 14\%$ of the nominal loading level $d_{0}$. By inspection, we see that $\rho_\infty$, when used with the PUS uncertainty set, performs much better at characterizing flexibility than when used on its corresponding box uncertainty set. The reason for this assessment is that $\rho_\infty$ spans a much wider range under PUS than under the box set. This is a desirable feature because system operators would have a more sensitive assessment of prevailing levels of flexibility. In the case of the box set, variations in $\rho_\infty$ are so narrow ($\approx 0.25$ about the base index, in comparison to $\approx 0.55$ with the PUS). Moreover, it does not even get close to one as the system loading becomes more critical, a feature required to signal clearly flexibility inadequacy.

\subsubsection{Level of Uncertainty Influence on Constraints Projections into the Demand Space}
We consider next the RTS nominal loading level with different levels of uncertainty $\eta = \{0.033, 0.067, 0.1\}$ for modelling load variations and a correlation factor $\alpha = 0.7$. The line capacities are reduced by half to have a more congested network. 
In Table.~\ref{table:05} we notice an increasing trend in umbrella line flow constraints when the level of uncertainty grows. This happens because the system has to cover wider ranges of residual demands which tends to increase the possibility of hitting line flow limit. 

In particular, for the case when $\eta = 0.067$ or $0.1$, the box set retains a higher percentage of its line constraints in comparison with the PUS, by 77.78\% and 66.67\%, respectively. On the other hand, in the generation-demand space most generator constraints are expected to be umbrellas. Still, while considering the highest level uncertainty, the PUS retains only up to 22.7\% of the total original constraints in comparison to 30.9\% with the conventional box approach. This is a positive outcome indicating the superiority of the PUS approach. With fewer remaining constraints, the evaluation of $\rho_r$ and RDC through DDIO is expected to be simpler and faster.    



\begin{table}
\captionsetup{font=footnotesize}
\caption{\sc{Percentage of Constraints Retained in the BA and its Corresponding Loadability Set After the Application of UCD}} 
\centering 
\scalebox{0.65}{
\begin{tabular}{ccccc} 
\hline
 & PUS & Box & PUS & Box \\ [1.0ex]
\hline 
$\eta$ & \multicolumn{2}{c}{Line constraints ($\%$ change)} & \multicolumn{2}{c}{Generation-demand space ($\%$ change)} \\ [1.0ex] 
\hline
0.033 & 9.21   & 9.21  (0.0\%)     & 20.6  & 23.7  (+15\%)   \\[1ex]
0.067 & 11.84  & 21.05 (+77.78\%)  & 20.6  & 27.8 (+34.95\%)   \\[1ex]
0.100 & 15.78  & 26.31 (+66.67\%)  & 22.7  & 30.9 (+36.12\%)  \\[1ex]
\hline
\end{tabular}}
\label{table:05} 
\end{table}

\section{Conclusion}\label{CONCLUSION}

This paper proposed a novel data-driven inverse optimization scheme for assessing flexibility explicitly in low-carbon power systems. Using historical demand data and its forecasts, polyhedral uncertainty sets can capture the spatial correlation of residual demands and their forecasting errors. A new flexibility metric provides a useful numerical and visualization tool to assess how residual demand forecast errors are expected to be handled given committed flexible resources. Moreover, umbrella constraint discovery is used to determine the minimum number of constraints that shape the feasibility region and how that can be affected under correlated uncertainties. The framework unlocks the integration assessment of renewable generation with loadability set approaches by defining critical constraints. Several interesting potential applications could be explored on the basis on the inverse optimization scheme presented in the paper. For example, one could assess how the flexibility metric $\rho_r$ and RDC are affected by deploying energy storage assets or demand response programs.


\bibliography{References} 
\bibliographystyle{ieeetr}
\end{document}